# CENTRAL LIMIT THEOREMS FOR SEQUENCES OF MULTIPLE STOCHASTIC INTEGRALS

By David Nualart and Giovanni Peccati

*Universitat de Barcelona and Université de Paris VI*

We characterize the convergence in distribution to a standard normal law for a sequence of multiple stochastic integrals of a fixed order with variance converging to 1. Some applications are given, in particular to study the limiting behavior of quadratic functionals of Gaussian processes.

**1. Introduction.** In this paper, we characterize the convergence in distribution to a normal $N(0,1)$ law for a sequence of random variables $F_k$ belonging to a fixed Wiener chaos, and with variance tending to 1. We show that a necessary and sufficient condition for this convergence is that the moment of fourth order of $F_k$ converges to 3. If $F_k$ is a multiple stochastic integral of order $n$ of a symmetric and square integrable kernel $f_k$, for instance on $[0,1]^n$, another necessary and sufficient condition for the above convergence is that, for all $p = 1, \ldots, n-1$, the contractions of order $p$ (defined by $f_k^{\otimes p} := f_k \otimes_p f_k$) converge to zero in $L^2([0,1]^{2(n-p)})$ as $k$ tends to infinity.

In general, we call "central limit theorem" (CLT in the sequel) any result describing the weak convergence of a sequence of nonlinear functionals of a Gaussian process (or of a Gaussian measure) toward a standard normal law. The reader is referred to Major (1981), Maruyama (1982, 1985), Giraitis and Surgailis (1985) and the references therein for results in this direction. Here, we shall observe that, in the above quoted references, the authors establish sufficient conditions to have a CLT in the general case of sequences of functionals having a possibly infinite Wiener–Itô expansion. A common technique used in such a study is the *method of moments* [see, e.g., Maruyama (1985)], requiring a determination of all moments associated to a given functional, usually estimated by means of the so-called *diagram formulae* [see Surgailis (2000) for a detailed survey]. On the other hand, our techniques (which are mainly based on a stochastic calculus result due to Dambis, Dubins and Schwarz [see Revuz and Yor (1999), Chapter V and Section 3.1]) naturally bring the need to estimate and control expressions related uniquely to the fourth moment of each element of the sequence $F_k$. To this end, we apply extensively some version of the *product formula* for multiple stochastic integrals, such as the one presented in Nualart [(1995), Proposition 1.1.3], and perform calculations that are very close in spirit to the ones contained in the first part of Üstünel and Zakai (1989).

Our results are specifically motivated by recent works on limit theorems for quadratic functionals of Brownian motion and Brownian bridge [see Deheuvels and Martynov (2004) and Peccati and Yor (2004a, b)], as well as Brownian sheet and related processes [see Deheuvels, Peccati and Yor (2004)]. We provide examples and applications, mainly related to quadratic functionals of a fractional Brownian motion, with Hurst parameter $H > \frac{1}{2}$, and of a standard Brownian sheet.

**2. The main result.** Consider a separable Hilbert space $H$. Let $\{e_k : k \geq 1\}$ be a complete orthonormal system in $H$. For every $n \geq 1$, we denote $H^{\odot n}$ the $n$th symmetric tensor product of $H$. For $p = 0, \ldots, n$, and for every $f \in H^{\odot n}$, we define the *contraction* of $f$ of order $p$ to be the element of $H^{\otimes 2(n-p)}$ defined by

$$f^{\otimes p} = \sum_{i_1, \ldots, i_p = 1}^{\infty} \langle f, e_{i_1} \otimes \cdots \otimes e_{i_p} \rangle_{H^{\otimes p}} \otimes \langle f, e_{i_1} \otimes \cdots \otimes e_{i_p} \rangle_{H^{\otimes p}},$$

and we denote by $(f^{\otimes p})_s$ its symmetrization.

In what follows, we will write

$$X = \{X(h) : h \in H\}$$

for an *isonormal Gaussian process* on $H$. This means that $X$ is a centered Gaussian family indexed by the elements of $H$, defined on some probability space $(\Omega, \mathcal{F}, \mathbb{P})$ and such that, for every $h, h' \in H$,

$$\mathbb{E}(X(h)X(h')) = \langle h, h' \rangle_H.$$

For every $n \geq 1$, we will denote by $I_n^X$ the isometry between $H^{\odot n}$ equipped with the norm $\sqrt{n!} \| \cdot \|_{H^{\otimes n}}$ and the $n$th Wiener chaos of $X$.

In the particular case where $H = L^2(A, \mathcal{A}, \mu)$, $(A, \mathcal{A})$ is a measurable space and $\mu$ is a $\sigma$-finite and nonatomic measure, then $H^{\odot n} = L_s^2(A^n, \mathcal{A}^{\otimes n}, \mu^{\otimes n})$ is the space of symmetric and square integrable functions on $A^n$ and for every $f \in H^{\odot n}$, $I_n^X(f)$ is the *multiple Wiener–Itô integral* (of order $n$) of $f$ with respect to $X$, as defined, for example, in Nualart [(1995), Section 1.1.2].

Our main result is the following.

THEOREM 1. *Let the above notation and assumptions prevail, and fix $n \geq 2$. Then, for any sequence of elements $\{f_k : k \geq 1\}$ such that $f_k \in H^{\odot n}$ for every $k$, and*

$$\lim_{k \to +\infty} n! \|f_k\|^2_{H^{\otimes n}} = \lim_{k \to +\infty} \mathbb{E}[I_n^X(f_k)^2] = 1, \tag{1}$$

*the following conditions are equivalent*:

(i) $\lim_{k \to +\infty} \mathbb{E}[I_n^X(f_k)^4] = 3$;
(ii) *for every* $p = 1, \ldots, n-1$, $\lim_{k \to +\infty} \|f_k^{\otimes p}\|^2_{H^{\otimes 2(n-p)}} = 0$;
(iii) *as $k$ goes to infinity, the sequence $\{I_n^X(f_k) : k \geq 1\}$ converges in distribution to a standard Gaussian random variable.*

As a consequence of this result we obtain

COROLLARY 2. *Fix $n \geq 2$ and $f \in H^{\odot n}$ such that $\mathbb{E}[I_n^X(f)^2] = 1$. Then the distribution of $I_n^X(f)$ cannot be normal and $\mathbb{E}[I_n^X(f)^4] \neq 3$.*

PROOF. If $I_n^X(f)$ had a normal distribution or $\mathbb{E}[I_n^X(f)^4] = 3$, then, according to Theorem 1, for every $p = 1, \ldots, n-1$, we would have $f^{\otimes p} = 0$. Thus, for each $v \in H^{\otimes(n-p)}$ we obtain

$$0 = \langle f^{\otimes p}, v \otimes v \rangle_{H^{\otimes 2(n-p)}} = \|\langle f, v \rangle_{H^{\otimes(n-p)}}\|^2_{H^{\otimes p}},$$

which implies $f = 0$. □

REMARK 1. The fact that a random variable with the form $I_n^X(f)$, $n > 1$, cannot be Gaussian is well known. See, for example, the discussion contained in Janson [(1997), Chapter VI].

To prove Theorem 1 in its most general form, we shall first deal with the case of $X$ being the Gaussian family generated by a standard Brownian motion on $[0, 1]$. To this end, some further notation is needed.

2.1. *The Brownian case.* For every $n \geq 1$ and $T > 0$, given a permutation $\pi$ of $(1, \ldots, n)$, we set

$$\Delta^n_{\pi, T} = \{(t_1, \ldots, t_n) \in \mathbb{R}^n : T > t_{\pi(1)} > \cdots > t_{\pi(n)} > 0\},$$

and we define

$$\Delta^n_{\pi_0, T} := \Delta^n_T = \{(t_1, \ldots, t_n) \in \mathbb{R}^n : T > t_1 > \cdots > t_n > 0\}$$

to be the simplex contained in $[0, T]^n$. Now, for a given $n \geq 1$, take $f \in L^2(\Delta^n_T, dt_1 \cdots dt_n) = L^2(\Delta^n_T)$. To such an $f$ we associate the symmetric function on $[0, T]^n$,

$$\tilde{f}(t_1, \ldots, t_n) = \sum_\pi f(t_{\pi(1)}, \ldots, t_{\pi(n)}) \mathbf{1}_{\Delta^n_{\pi, T}}(t_1, \ldots, t_n),$$

where $\pi$ runs over all permutations of $(1,\ldots,n)$. For $p=0,\ldots,n$, the contraction of $\tilde{f}$ of order $p$ on $[0,T]$ is the application

$$(t_1,\ldots,t_{2n-2p}) \mapsto \tilde{f}^{\otimes p,T}(t_1,\ldots,t_{2n-2p})$$
$$= \int_{[0,T]^p} \tilde{f}(u_1,\ldots,u_p,t_1,\ldots,t_{n-p})$$
$$\times \tilde{f}(u_1,\ldots,u_p,t_{n-p+1},\ldots,t_{2n-2p})\,du_1\cdots du_p$$

and we note $\tilde{f}^{\otimes p,T}(\cdot)_s$ the symmetrization of $\tilde{f}^{\otimes p,T}$.

In this section, we note $W = \{W_t : \in [0,1]\}$ a standard Brownian motion initialized at zero. Note that, in the terminology of the previous paragraph, the centered Gaussian space generated by $W$ can be identified with an isonormal Gaussian process on $H = L^2([0,1], dt)$. For every $n$, for every $f \in L^2(\Delta_1^n)$ and every $t \in (0,1]$, we put

$$J_n^t(f) = \int_0^t dW_{s_1} \cdots \int_0^{s_{n-1}} dW_{s_n} f(s_1,\ldots,s_n)$$

and also $I_n^t(\tilde{f}) = n! J_n^t(f)$, so that

$$I_n^1(\tilde{f}) = I_n^X(\tilde{f}),$$

where $X$ is once again the isonormal process generated by $W$. The following result translates the content of Theorem 1 in the context of this section.

PROPOSITION 3.  *Let the above notation prevail, fix $n \geq 2$ and consider a collection $\{g_k : k \geq 1\}$ of elements of $L^2(\Delta_1^n)$ such that*

$$\lim_{k \to +\infty} \|g_k\|_{L^2(\Delta_1^n)} := \lim_{k \to +\infty} \|g_k\|_{\Delta_1^n} = 1. \tag{2}$$

*Then, the following conditions are equivalent:*

(i) $\lim_{k \to +\infty} \mathbb{E}[J_n^1(g_k)^4] = 3$;

(ii) *for every* $p=1,\ldots,n-1$, $\lim_{k \to +\infty} \|\tilde{g}_k^{\otimes p,1}\|_{[0,1]^{2n-2p}}^2 = 0$;

(iii) *as* $k \to \infty$, *the sequence* $\{J_n^1(g_k) : k \geq 1\}$ *converges in distribution to a standard Gaussian random variable;*

(iv) *for every* $p = 0,\ldots,n-2$,

$$\lim_{k \to \infty} \int_{\Delta_1^{2(n-1-p)}} ds_1 \cdots ds_{2(n-1-p)} \left( \int_{s_1}^1 dt\, \tilde{g}_{k,t}^{\otimes p, t}(s_1,\ldots,s_{2(n-1-p)})_s \right)^2 = 0,$$
(3)
*where, for every fixed $t \in [0,1]$, $g_{k,t}$ stands for the function on $\Delta_t^{n-1}$ given by*

$$(s_1,\ldots,s_{n-1}) \mapsto g_k(t,s_1,\ldots,s_{n-1}).$$

PROOF. We will prove the following implications:

$$(iv) \Rightarrow (iii) \Rightarrow (i) \Rightarrow (ii) \Rightarrow (iv).$$

(iv) $\Rightarrow$ (iii). Suppose that (3) holds. According to the Dambis–Dubins–Schwarz theorem [see Revuz and Yor (1999), Chapter V], for every $k$ there exists a standard Brownian motion $W^{(k)}$ such that

$$(4) \qquad J_n^1(g_k) = W^{(k)}_{\int_0^1 dt (J_{n-1}^t(g_{k,t}))^2} = W^{(k)}_{[(n-1)!]^{-2} \int_0^1 dt (I_{n-1}^t(\tilde{g}_{k,t}))^2}.$$

We can expand the stochastic time-change in (4) by means of a multiplication formula for multiple Wiener integrals [see, e.g., Nualart (1995), Proposition 1.1.3] to obtain that

$$\int_0^1 dt \, (I_{n-1}^t(\tilde{g}_{k,t}))^2 = \int_0^1 dt \left[ \sum_{p=0}^{n-1} c_{n-1,p} I_{2(n-1-p)}^t(\tilde{g}_{k,t}^{\otimes_p, t}(\cdot)_s) \right],$$

where $c_{n-1,p} = p! \binom{n-1}{p}^2$, and therefore

$$(5) \qquad \frac{\int_0^1 dt (I_{n-1}^t(\tilde{g}_{k,t}))^2}{[(n-1)!]^2} = \|g_k\|_{\Delta_1^n}^2 + \frac{\sum_{p=0}^{n-2} c_{n-1,p} \int_0^1 dt (I_{2(n-1-p)}^t(\tilde{g}_{k,t}^{\otimes_p, t}(\cdot)_s))}{[(n-1)!]^2}.$$

Moreover, thanks to a stochastic Fubini theorem for multiple integrals, we obtain that for any $p = 0, \ldots, n-2$,

$$\frac{\int_0^1 dt (I_{2(n-1-p)}^t(\tilde{g}_{k,t}^{\otimes_p, t}(\cdot)_s))}{[2(n-1-p)]!}$$

$$= \int_0^1 dt (J_{2(n-1-p)}^t(\tilde{g}_{k,t}^{\otimes_p, t}(\cdot)_s))$$

$$= \int_{\Delta_1^{2(n-1-p)}} dW_{s_1} \cdots dW_{s_{2(n-1-p)}} \left( \int_{s_1}^1 dt \, \tilde{g}_{k,t}^{\otimes_p, t}(s_1, \ldots, s_{2(n-1-p)})_s \right).$$

Now take the $L^2$-norm of the right-hand side of (5), let $k$ go to infinity and observe that if (2) and (3) are verified, then the pair

$$\left( W_\cdot^{(k)}, [(n-1)!]^{-2} \int_0^1 dt \, (I_{n-1}^t(\tilde{g}_{k,t}))^2 \right)$$

converges weakly to $(W_\cdot, 1)$, so that the conclusion is immediately achieved by using, for example, Theorem 3.1 in Whitt (1980), as well as formula (4).

(iii) $\Rightarrow$ (i). Trivial, given condition (2).

(i) $\Rightarrow$ (ii). Observe first that, for every $k$,

$$\mathbb{E}[J_n^1(g_k)^4] = (n!)^{-4} \mathbb{E}[I_n^1(\tilde{g}_k)^4].$$

Now, as

$$I_n^1(\tilde{g}_k)^2 = \sum_{p=0}^n p! \binom{n}{p}^2 I_{2n-2p}^1(\tilde{g}_k^{\otimes p,1}(\cdot)_s)$$

$$= n!\|\tilde{g}_k\|_{[0,1]^n}^2 + I_{2n}^1(\tilde{g}_k^{\otimes 0,1}(\cdot)_s) + \sum_{p=1}^{n-1} p! \binom{n}{p}^2 I_{2n-2p}^1(\tilde{g}_k^{\otimes p,1}(\cdot)_s),$$

due again to the multiplication formula, we obtain immediately

$$\mathbb{E}[J_n^1(g_k)^4] = \|g_k\|_{\Delta_n^1}^4 + \frac{(2n)!}{(n!)^4}\|\tilde{g}_k^{\otimes 0,1}(\cdot)_s\|_{[0,1]^{2n}}^2$$

$$+ (n!)^{-4} \sum_{p=1}^{n-1} \left[p!\binom{n}{p}^2\right]^2 (2n-2p)!\|\tilde{g}_k^{\otimes p,1}(\cdot)_s\|_{[0,1]^{2n-2p}}^2.$$

Now call $\Pi_{2n}$ the set of the $(2n)!$ permutations of the set $(1,\ldots,2n)$, whose generic elements will be denoted by $\pi$, $\pi'$, and so on. On such a set, we introduce the following notation: for $p=0,\ldots,n$, we write

$$\pi \overset{p}{\sim} \pi'$$

if the set $(\pi(1),\ldots,\pi(n)) \cap (\pi'(1),\ldots,\pi'(n))$ contains exactly $p$ elements. Note that, for a given $\pi \in \Pi_{2n}$, there are $\binom{n}{p}^2 (n!)^2$ permutations $\pi'$ such that

$$\pi' \overset{p}{\sim} \pi.$$

Moreover, it is easily seen that if $\pi' \overset{0}{\sim} \pi$ or $\pi' \overset{n}{\sim} \pi$, then

$$\int_{[0,1]^{2n}} da_1 \cdots da_{2n} \tilde{g}_k^{\otimes 0,1}(a_{\pi(1)},\ldots,a_{\pi(2n)}) \tilde{g}_k^{\otimes 0,1}(a_{\pi'(1)},\ldots,a_{\pi'(2n)})$$

$$= \|\tilde{g}_k\|_{[0,1]^n}^4 = (n!)^2 \|g_k\|_{\Delta_n^1}^4,$$

so that

$$\frac{(2n)!}{(n!)^4}\|\tilde{g}_k^{\otimes 0,1}(\cdot)_s\|_{[0,1]^{2n}}^2$$

$$= \frac{1}{(n!)^4 (2n)!} \sum_{\pi \in \Pi_{2n}} \left[\|g_k\|_{\Delta_n^1}^4 \left(\sum_{\pi' \overset{0}{\sim} \pi} (n!)^2 + \sum_{\pi' \overset{n}{\sim} \pi} (n!)^2\right)\right.$$

$$\left.+ \sum_{p=1}^{n-1} \sum_{\pi' \overset{p}{\sim} \pi} \int_{[0,1]^{2n}} da_1 \cdots da_{2n} \tilde{g}_k^{\otimes 0,1}(a_\pi) \tilde{g}_k^{\otimes 0,1}(a_{\pi'})\right],$$

where $a_\pi = (a_{\pi(1)},\ldots,a_{\pi(2n)})$. Since, for $p=1,\ldots,n-1$,

$$\sum_{\pi' \overset{p}{\sim} \pi} \int_{[0,1]^{2n}} da_1 \cdots da_{2n} \tilde{g}_k^{\otimes 0,1}(a_\pi) \tilde{g}_k^{\otimes 0,1}(a_{\pi'}) = \binom{n}{p}^2 (n!)^2 \|\tilde{g}_k^{\otimes p,1}\|_{[0,1]^{2n-2p}}^2,$$

we get

$$\frac{(2n)!}{(n!)^4}\|\tilde{g}_k^{\otimes_0,1}(\cdot)_s\|_{[0,1]^{2n}}^2 = 2\|g_k\|_{\Delta_n^1}^4 + \sum_{p=1}^{n-1}\frac{\|\tilde{g}_k^{\otimes_p,1}\|_{[0,1]^{2n-2p}}^2}{(p!(n-p)!)^2},$$

and therefore

$$\mathbb{E}[J_n^1(g_k)^4] = 3\|g_k\|_{\Delta_n^1}^4$$
$$+ \sum_{p=1}^{n-1}(p!(n-p)!)^{-2}$$
$$\times \left[\|\tilde{g}_k^{\otimes_p,1}\|_{[0,1]^{2n-2p}}^2 + \binom{2n-2p}{n-p}\|\tilde{g}_k^{\otimes_p,1}(\cdot)_s\|_{[0,1]^{2n-2p}}^2\right].$$

This yields in particular that, if $\mathbb{E}[J_n^1(g_k)^4]$ converges to 3, and (2) holds, then necessarily, for every $p = 1, \ldots, n-1$,

(6) $$\lim_k \|\tilde{g}_k^{\otimes_p,1}\|_{[0,1]^{2n-2p}}^2 = 0.$$

(ii) $\Rightarrow$ (iv). We introduce some notation: for any $m \geq 1$, $\mathbf{x}_m$ is shorthand for a vector $(x_1, \ldots, x_m) \in \mathbb{R}^m$, $\hat{\mathbf{x}}_m = \max_i(x_i)$, and $d\mathbf{x}_m$ stands for Lebesgue measure on $\mathbb{R}^m$. To conclude the proof, we shall show that, for $p = 1, \ldots, n-1$, condition (6) implies necessarily that

$$\lim_k \int_{[0,1]^{n-p}} d\mathbf{s}_{n-p} \int_{[0,1]^{n-p}} d\boldsymbol{\tau}_{n-p} \left(\int_{\hat{\mathbf{s}}_{n-p}\vee\hat{\boldsymbol{\tau}}_{n-p}}^1 dt\, \tilde{g}_{k,t}^{\otimes_{p-1},t}(\mathbf{s}_{n-p}, \boldsymbol{\tau}_{n-p})\right)^2 = 0.$$

Now, since for every $(t, s_2, \ldots, s_n) \in [0,1]^n$,

$$\tilde{g}_{k,t}(s_2, \ldots, s_n) = \mathbf{1}_{[0,t]^{n-1}}(s_2, \ldots, s_n)\tilde{g}_k(t, s_2, \ldots, s_n),$$

we obtain that

$$\int_{[0,1]^{n-p}} d\mathbf{s}_{n-p} \int_{[0,1]^{n-p}} d\boldsymbol{\tau}_{n-p} \left(\int_{\hat{\mathbf{s}}_{n-p}\vee\hat{\boldsymbol{\tau}}_{n-p}}^1 dt\, \tilde{g}_{k,t}^{\otimes_{p-1},t}(\mathbf{s}_{n-p}, \boldsymbol{\tau}_{n-p})\right)^2$$
$$= \int_0^1 dt' \int_{[0,t']^{p-1}} d\mathbf{v}_{p-1} \int_0^1 dt \int_{[0,t]^{p-1}} d\mathbf{u}_{p-1} \mathbf{1}_{(\hat{\mathbf{u}}_{p-1} \leq t, \hat{\mathbf{v}}_{p-1} \leq t')}$$
$$\times \left(\int_{[0,t\wedge t']^{n-p}} d\mathbf{s}_{n-p}\, \tilde{g}_k(t, \mathbf{u}_{p-1}, \mathbf{s}_{n-p})\tilde{g}_k(t', \mathbf{v}_{p-1}, \mathbf{s}_{n-p})\right)^2$$
$$\leq \int_0^1 dt' \int_{[0,1]^{p-1}} d\mathbf{v}_{p-1} \int_0^1 dt \int_{[0,1]^{p-1}} d\mathbf{u}_{p-1}$$
$$\times \left(\int_{[0,t\wedge t']^{n-p}} d\mathbf{s}_{n-p}\, \tilde{g}_k(t, \mathbf{u}_{p-1}, \mathbf{s}_{n-p})\tilde{g}_k(t', \mathbf{v}_{p-1}, \mathbf{s}_{n-p})\right)^2$$

$$= \int_{[0,1]^{n-p}} d\mathbf{s}_{n-p} \int_{[0,1]^{n-p}} d\boldsymbol{\tau}_{n-p}$$

$$\times \left( \int_{\hat{\mathbf{s}}_{n-p} \vee \hat{\boldsymbol{\tau}}_{n-p}}^{1} dt \int_{[0,1]^{p-1}} d\mathbf{u}_{p-1}\, \tilde{g}_k(t, \mathbf{u}_{p-1}, \mathbf{s}_{n-p}) \tilde{g}_k(t, \mathbf{u}_{p-1}, \boldsymbol{\tau}_{n-p}) \right)^2.$$

Now, since condition (6) holds by assumption, we know that

$$A(k) = \int_{[0,1]^{n-p}} d\mathbf{s}_{n-p} \int_{[0,1]^{n-p}} d\boldsymbol{\tau}_{n-p} \left( \int_{[0,1]^p} d\mathbf{u}_p\, \tilde{g}_k(\mathbf{u}_p, \mathbf{s}_{n-p}) \tilde{g}_k(\mathbf{u}_p, \boldsymbol{\tau}_{n-p}) \right)^2$$

converges to zero, as $k$ goes to infinity. We may expand the above expression and write

$$A(k) = \int_{[0,1]^{n-p}} d\mathbf{s}_{n-p} \int_{[0,1]^{n-p}} d\boldsymbol{\tau}_{n-p}$$

$$\times \left( \int_0^{\hat{\mathbf{s}}_{n-p} \vee \hat{\boldsymbol{\tau}}_{n-p}} dt \int_{[0,1]^{p-1}} d\mathbf{u}_{p-1}\, \tilde{g}_k(t, \mathbf{u}_{p-1}, \mathbf{s}_{n-p}) \tilde{g}_k(t, \mathbf{u}_{p-1}, \boldsymbol{\tau}_{n-p}) \right)^2$$

$$+ \int_{[0,1]^{n-p}} d\mathbf{s}_{n-p} \int_{[0,1]^{n-p}} d\boldsymbol{\tau}_{n-p}$$

$$\times \left( \int_{\hat{\mathbf{s}}_{n-p} \vee \hat{\boldsymbol{\tau}}_{n-p}}^{1} dt \int_{[0,1]^{p-1}} d\mathbf{u}_{p-1}\, \tilde{g}_k(t, \mathbf{u}_{p-1}, \mathbf{s}_{n-p}) \tilde{g}_k(t, \mathbf{u}_{p-1}, \boldsymbol{\tau}_{n-p}) \right)^2$$

$$+ 2 \int_{[0,1]^{n-p}} d\mathbf{s}_{n-p} \int_{[0,1]^{n-p}} d\boldsymbol{\tau}_{n-p}$$

$$\times \left( \int_0^{\hat{\mathbf{s}}_{n-p} \vee \hat{\boldsymbol{\tau}}_{n-p}} dt \int_{[0,1]^{p-1}} d\mathbf{u}_{p-1}\, \tilde{g}_k(t, \mathbf{u}_{p-1}, \mathbf{s}_{n-p}) \tilde{g}_k(t, \mathbf{u}_{p-1}, \boldsymbol{\tau}_{n-p}) \right)$$

$$\times \left( \int_{\hat{\mathbf{s}}_{n-p} \vee \hat{\boldsymbol{\tau}}_{n-p}}^{1} dt' \int_{[0,1]^{p-1}} d\mathbf{v}_{p-1}\, \tilde{g}_k(t', \mathbf{v}_{p-1}, \mathbf{s}_{n-p}) \tilde{g}_k(t', \mathbf{v}_{p-1}, \boldsymbol{\tau}_{n-p}) \right)$$

$$= A_1(k) + A_2(k) + A_3(k).$$

In particular, the Fubini theorem yields

$$A_3(k) = \int_{[0,1]^p} dt\, d\mathbf{u}_{p-1} \int_{[0,1]^p} dt'\, d\mathbf{v}_{p-1}$$

$$\times \left( \int_{[0,1]^{n-p}} d\mathbf{s}_{n-p}\, \mathbf{1}_{(t' \wedge t \leq \hat{\mathbf{s}}_{n-p} \leq t' \vee t)}\, \tilde{g}_k(t, \mathbf{u}_{p-1}, \mathbf{s}_{n-p}) \tilde{g}_k(t', \mathbf{v}_{p-1}, \mathbf{s}_{n-p}) \right)^2$$

$$+ 2 \int_{[0,1]^p} dt\, d\mathbf{u}_{p-1} \int_{[0,1]^p} dt'\, d\mathbf{v}_{p-1}$$

$$\times \Bigg( \int_{[0,1]^{n-p}} d\mathbf{s}_{n-p}\, \mathbf{1}_{(t' \wedge t \leq \hat{\mathbf{s}}_{n-p} \leq t' \vee t)}$$

$$\times \tilde{g}_k(t, \mathbf{u}_{p-1}, \mathbf{s}_{n-p}) \tilde{g}_k(t', \mathbf{v}_{p-1}, \mathbf{s}_{n-p})\Big)$$

$$\times \left( \int_{[0,1]^{n-p}} d\boldsymbol{\tau}_{n-p} \mathbf{1}_{(\hat{\boldsymbol{\tau}}_{n-p} \leq t' \wedge t)} \tilde{g}_k(t, \mathbf{u}_{p-1}, \boldsymbol{\tau}_{n-p}) \tilde{g}_k(t', \mathbf{v}_{p-1}, \boldsymbol{\tau}_{n-p}) \right).$$

Since

$$A_2(k) = \int_{[0,1]^p} dt\, d\mathbf{u}_{p-1} \int_{[0,1]^p} dt'\, d\mathbf{v}_{p-1}$$

$$\times \left( \int_{[0,1]^{n-p}} d\boldsymbol{\tau}_{n-p} \mathbf{1}_{(\hat{\boldsymbol{\tau}}_{n-p} \leq t' \wedge t)} \tilde{g}_k(t, \mathbf{u}_{p-1}, \boldsymbol{\tau}_{n-p}) \tilde{g}_k(t', \mathbf{v}_{p-1}, \boldsymbol{\tau}_{n-p}) \right)^2,$$

we obtain that, for every $k$, $A_2(k) + A_3(k) \geq 0$, and therefore that $A_1(k)$ must converge to zero as $k$ tends to infinity. This gives immediately the desired conclusion. □

REMARK 2. There exists an elegant explanation of the implication (iii) ⇒ (iv) for the case $n = 2$. To do this, consider for simplicity a family $\{I_2^1(g_k)\}_{k \geq 1}$ of multiple integrals with symmetric kernels $\{g_k\}_{k \geq 1}$ such that $\|g_k\|_{[0,1]^2} = 1$. Then, by using the formula

$$\mathbb{E}[\exp(i\lambda I_2^1(g_k))] = \exp\left[ \sum_{j \geq 2} (-1)^{j+1} \frac{(-2i\lambda)^j}{j} \operatorname{Tr}(g_k^j) \right],$$

where $g_k^j$ is the $j$th power of the Hilbert–Schmidt operator associated to the kernel $g_k(\cdot, \cdot)$, it is easy to show that if $I_2^1(g_k)$ converges in law to a standard Gaussian random variable as $k$ goes to infinity, then necessarily

$$\lim_k \operatorname{Tr}(g_k^j) = 0$$

for every $j > 2$. In particular, by taking $j = 4$, we obtain that

$$\lim_k \operatorname{Tr}(g_k^4) = \lim_k \int_0^1 dt \int_0^1 du \left( \int_0^1 g_k(t, r) g_k(u, r)\, dr \right)^2 = 0.$$

But some calculations (analogous to the ones at the end of the proof of Proposition 3) yield

$$\int_0^1 dt \int_0^1 du \left( \int_0^1 g_k(t, r) g_k(u, r)\, dr \right)^2$$

$$= \int_0^1 dt \int_0^1 du \left( \int_0^{t \vee u} g_k(t, r) g_k(u, r)\, dr \right)^2$$

$$+ \int_0^1 dr \int_0^1 dr' \left( \int_0^{r \wedge r'} dt\, g_k(t, r) g_k(t, r') \right)^2$$

$$+ \int_0^1 dr \int_0^1 dr' \left( \int_{r \wedge r'}^{r \vee r'} dt\, g_k(t,r) g_k(t,r') \right)^2$$

$$+ 2 \int_0^1 dr \int_0^1 dr' \left( \int_{r \wedge r'}^{r \vee r'} dt\, g_k(t,r) g_k(t,r') \right) \left( \int_0^{r \wedge r'} dt\, g_k(t,r) g_k(t,r') \right),$$

thus giving the desired conclusion.

A slight refinement of Proposition 3 is the following.

PROPOSITION 4. *Let $\{J_n^1(g_k) : k \geq 1\}$ be a sequence of iterated integrals as in the statement of Proposition* 3 *satisfying either condition* (i), (ii), (iii) *or* (iv). *Suppose, moreover, that the sequence*

$$\left\{ \int_0^1 dt\, |I_{n-1}^t(\tilde{g}_{k,t})| : k \geq 1 \right\}$$

*converges to zero in probability as $k$ goes to infinity. Then, the pair*

$$(J_n^1(g_k), W)$$

*converges in distribution to*

$$(N, W),$$

*where $N$ is a standard Gaussian random variable independent of $W$.*

PROOF. We shall only prove the asymptotic independence, which in this case is given by an application of an asymptotic version of Knight's theorem, such as the one stated, for instance, in Revuz and Yor [(1999), Chapter XIII]. □

2.2. *Proof of Theorem* 1. Let the assumptions and notation of Theorem 1 prevail. Since $H$ is a separable Hilbert space, for every $n \geq 1$ there exists an application $i_n(\cdot)$ from $H^{\odot n}$ onto $L_s^2([0,1]^n, dt_1 \cdots dt_n)$, such that, for every $n$, $i_n(\cdot)$ is an isometry and, moreover, the following equality holds for every $f \in H^{\odot n}$:

$$n! J_n^1(i_n(f)) = n! \int_0^1 \cdots \int_0^{t_{n-1}} dW_{t_1} \cdots dW_{t_n}\, i_n(f)(t_1, \ldots, t_n) \stackrel{\text{(law)}}{=} I_n^X(f),$$

where $W$ is a standard Brownian motion on $[0,1]$. It is therefore clear that the sequence $\{f_k : k \geq 1\}$ in the statement of Theorem 1 is such that:

1. $\lim_{k \to +\infty} n! \|i_n(f_k)\|^2_{[0,1]^n} = 1$.
2. $\lim_k \mathbb{E}[I_n^X(f_k)^4] = 3$ if, and only if,

$$(n!)^4 \lim_k \mathbb{E}[J_n^1(i_n(f_k))^4] = 3.$$

3. For every $p = 0, \ldots, n-1$, $\lim_k \|f_k^{\otimes p}\|_{H^{\odot 2n-2p}}^2 = 0$ if, and only if,
$$\lim_k \|i_n(f_k)^{\otimes p}\|_{[0,1]^{2n-2p}}^2 = 0.$$

4. For $k \to +\infty$, the convergence
$$I_n^X(f_k) \overset{\text{(law)}}{\Longrightarrow} N,$$
where $N$ is a standard Gaussian random variable, takes place if, and only if,
$$n! J_n^1(i_n(f_k)) \overset{\text{(law)}}{\Longrightarrow} N.$$

The proof is easily concluded.

## 3. Examples and applications.

3.1. *Quadratic functionals of a fractional Brownian motion.* Consider a fractional Brownian motion (fBm) $B^H = \{B_t^H, t \in [0,1]\}$ of Hurst parameter $H \in (0,1)$. That is, $B^H$ is a centered Gaussian process with the covariance function

(7) $$R_H(t,s) = \tfrac{1}{2}(s^{2H} + t^{2H} - |t-s|^{2H}).$$

According to the so-called *Jeulin's lemma* [see Jeulin (1980), Lemma 1, page 44], with probability 1, and for every $H$,

(8) $$\int_0^1 (B_t^H)^2 \frac{dt}{t^{2H+1}} = +\infty.$$

In what follows, we shall use the results of the previous sections to characterize the speed at which the two quantities

(9) $$F_\beta = \int_0^1 t^{2\beta}(B_t^H)^2\, dt \quad \text{and} \quad L_\varepsilon = \int_\varepsilon^1 (B_t^H)^2 \frac{dt}{t^{2H+1}}$$

diverge to infinity, respectively, when $\beta$ tends to $-H - \tfrac{1}{2}$, and when $\varepsilon$ tends to 0. In particular, to be able to apply Theorem 1, we will mainly concentrate on the case $H > \tfrac{1}{2}$.

We will first find the Wiener chaos expansion of $F_\beta$. Clearly,
$$\mathbb{E}(F_\beta) = \frac{1}{2\beta + 2H + 1}.$$

If $D$ denotes the derivative operator [see Alòs and Nualart (2003)], we have
$$D_s F_\beta = 2\int_s^1 t^{2\beta} B_t^H\, dt$$

and
$$D_r D_s F_\beta = \frac{2}{2\beta+1}(1-(s\vee r)^{2\beta+1}).$$

Stroock's formula [see Stroock (1987)] yields

$$(10) \qquad F_\beta = \frac{1}{2\beta+2H+1} + \frac{1}{2\beta+1} I_2(1-(\cdot\vee\cdot)^{2\beta+1}),$$

where $I_2$ denotes the double stochastic integral with respect to $B^H$.

The following result is another CLT.

PROPOSITION 5. *If $H > \frac{1}{2}$, then*

$$(11) \qquad (2\beta+2H+1)\int_0^1 t^{2\beta}(B_t^H)^2\,dt \xrightarrow[\beta\downarrow -H-1/2]{L^2(\Omega)} 1,$$

*and*

$$(12) \qquad (2\beta+2H+1)^{-1/2}\left((2\beta+2H+1)\int_0^1 t^{2\beta}(B_t^H)^2\,dt - 1\right) \xrightarrow[\beta\downarrow -H-1/2]{Law} c_H N,$$

*where $N$ is a random variable with distribution $N(0,1)$ and $c_H$ is a constant depending on $H$.*

PROOF. We denote by $\mathcal{E}$ the set of step functions on $[0,1]$. Let $\mathcal{H}$ be the Hilbert space defined as the closure of $\mathcal{E}$ with respect to the scalar product
$$\langle \mathbf{1}_{[0,t]}, \mathbf{1}_{[0,s]}\rangle_{\mathcal{H}} = R_H(t,s) = \alpha_H \int_0^t\int_0^s |r-u|^{2H-2}\,du\,dr,$$
where $\alpha_H = H(2H-1)$. The mapping $\mathbf{1}_{[0,t]} \to B_t$ can be extended to an isometry between $\mathcal{H}$ and the first chaos of $B^H$.

Let us compute
$$\mathbb{E}(I_2(1-(\cdot\vee\cdot)^{2\beta+1}))^2 = 2\|1-(\cdot\vee\cdot)^{2\beta+1}\|^2_{\mathcal{H}\otimes\mathcal{H}}$$
$$= 2\alpha_H^2 \int_{[0,1]^4} (1-(s\vee r)^{2\beta+1})(1-(t\vee u)^{2\beta+1})$$
$$\times |t-s|^{2H-2}|r-u|^{2H-2}\,ds\,dt\,dr\,du.$$

Simple computations lead to
$$\lim_{\beta\downarrow -H-1/2}(2\beta+2H+1)^2 \mathbb{E}(I_2(1-(\cdot\vee\cdot)^{2\beta+1}))^2 = 0,$$

because $\mathbb{E}(I_2(1-(s\vee r)^{2\beta+1}))^2$ behaves as $(2\beta+2H+1)^{-1}$ as $\beta\downarrow -H-\frac{1}{2}$. Therefore, (11) follows. In order to show (12), set
$$G_\beta = (2\beta+2H+1)^{-1/2}((2\beta+2H+1)F_\beta - 1)$$
$$= \frac{(2\beta+2H+1)^{1/2}}{2\beta+1} I_2(1-(\cdot\vee\cdot)^{2\beta+1}).$$

We have

$$\lim_{\beta \downarrow -H-1/2} \mathbb{E}(G_\beta^2)$$

$$= \lim_{\beta \downarrow -H-1/2} 2 \frac{(2\beta + 2H + 1)\alpha_H^2}{(2\beta + 1)^2}$$

$$\times \int_{[0,1]^4} (1 - (s \vee r)^{2\beta+1})(1 - (t \vee u)^{2\beta+1})$$

$$\times |t - s|^{2H-2} |r - u|^{2H-2} \, ds \, dt \, dr \, du$$

$$= \frac{\alpha_H^2}{2H^2} \lim_{\beta \downarrow -H-1/2} (2\beta + 2H + 1)$$

$$\times \int_{[0,1]^4} (s \vee r)^{2\beta+1} (t \vee u)^{2\beta+1}$$

$$\times |t - s|^{2H-2} |r - u|^{2H-2} \, ds \, dt \, dr \, du$$

$$= c_H^2.$$

Taking into account Theorem 1, it suffices to show that condition (ii) holds, that is,

$$\lim_{\beta \downarrow -H-1/2} \frac{(2\beta + 2H + 1)^2}{(2\beta + 1)^4}$$

$$\times \left\| \int_{[0,1]^2} (1 - (\cdot \vee r)^{2\beta+1}) \right.$$

$$\left. \times (1 - (\cdot \vee r')^{2\beta+1}) |r - r'|^{2H-2} \, dr \, dr' \right\|_{\mathcal{H} \otimes \mathcal{H}}^2 = 0.$$

The above norm can be written as

$$\alpha_H^4 \int_{[0,1]^8} (1 - (s \vee r)^{2\beta+1})(1 - (t \vee r')^{2\beta+1})$$

$$\times (1 - (s' \vee u)^{2\beta+1})(1 - (t' \vee u')^{2\beta+1})$$

$$\times |r - r'|^{2H-2} |u - u'|^{2H-2}$$

$$\times |s - s'|^{2H-2} |t - t'|^{2H-2} \, dr \, dr' \, du \, du' \, ds \, ds' \, dt \, dt',$$

and it is of the order of $(2\beta + 2H + 1)^{-1}$ as $\beta \downarrow -H - \frac{1}{2}$, which implies the desired result. $\square$

Note that we have also the following noncentral limit theorem.

PROPOSITION 6. *For any $H \in (0,1)$,*

$$(13) \qquad (2\beta + 2H + 1) \int_0^1 t^{2\beta} (B_t^H)^2 \, dt \xrightarrow[\beta \uparrow +\infty]{L^2(\Omega)} (B_1^H)^2.$$

PROOF. The convergence (13) follows from (10) and the fact that

$$I_2(1) = (B_1^H)^2 - 1. \qquad \square$$

REMARK 3. The combination of Propositions 5 and 6 generalizes results previously obtained for a standard Brownian motion and a standard Brownian bridge. In particular, analogs of formulae (11) and (13) for the case $H = \frac{1}{2}$ are proved in Deheuvels and Martynov (2004) and Peccati and Yor (2004a), whereas a version of (12) for $H = \frac{1}{2}$ is obtained in Peccati and Yor (2004a).

Next, simple calculations yield that the Wiener chaos decomposition of $L_\varepsilon$, for any $\varepsilon > 0$, is given by

$$L_\varepsilon = \log \frac{1}{\varepsilon} - \frac{1}{2H} I_2(1 - (\varepsilon \vee \cdot \vee \cdot)^{-2H}),$$

where the second summand on the right-hand side is the double stochastic integral, with respect to $B^H$, of the symmetric function

$$(s,t) \mapsto 1 - (\varepsilon \vee t \vee s)^{-2H}.$$

This yields the following CLT.

PROPOSITION 7. *If $H > \frac{1}{2}$, then*

$$(14) \qquad \left(\log \frac{1}{\varepsilon}\right)^{-1/2} \left(L_\varepsilon - \log \frac{1}{\varepsilon}\right) \xrightarrow[\varepsilon \downarrow 0]{Law} k_H N,$$

*where $N$ is a standard Gaussian random variable, and $k_H$ is a real constant depending exclusively on $H$.*

PROOF. We keep the notation in the proof of Proposition 5. Since for every $t > 0$ the Wiener chaos decomposition of $(B_t^H)^2$ is given by

$$(B_t^H)^2 = t^{2H} + I_2(\mathbf{1}_{(\cdot \vee \cdot \leq t)}),$$

it is easily verified that

$$\mathbb{E}[L_\varepsilon^2] = \int_\varepsilon^1 \frac{dt}{t^{2H+1}} \int_\varepsilon^1 \frac{ds}{s^{2H+1}} \mathbb{E}[(B_t^H)^2 (B_s^H)^2]$$

$$= \left(\log \frac{1}{\varepsilon}\right)^2 + \frac{1}{2} \int_\varepsilon^1 \frac{dt}{t^{2H+1}} \int_\varepsilon^1 \frac{ds}{s^{2H+1}} (t^{2H} + s^{2H} - |t-s|^{2H})^2$$

and consequently that

$$\mathbb{E}[L_\varepsilon^2] - \mathbb{E}[L_\varepsilon]^2 = \mathbb{E}[L_\varepsilon^2] - \left(\log\frac{1}{\varepsilon}\right)^2$$

behaves like $\log(1/\varepsilon)$ when $\varepsilon$ tends to zero. Thanks to Theorem 1, to prove (14) it is now sufficient to show that

$$\lim_{\varepsilon \downarrow 0}[\log(1/\varepsilon)]^{-2}$$
$$\times \left\|\int_{[0,1]^2}(1 - (\cdot \vee r \vee \varepsilon)^{-2H})(1 - (\cdot \vee r' \vee \varepsilon)^{-2H})|r-r'|^{2H-2}\,dr\,dr'\right\|^2_{\mathcal{H}\otimes\mathcal{H}}$$
$$= 0.$$

Now, it is clear that

$$\left\|\int_{[0,1]^2}(1 - (\cdot \vee r \vee \varepsilon)^{-2H})(1 - (\cdot \vee r' \vee \varepsilon)^{-2H})|r-r'|^{2H-2}\,dr\,dr'\right\|^2_{\mathcal{H}\otimes\mathcal{H}}$$
$$= \alpha_H^4 \int_{[0,1]^8}(1 - (s \vee r \vee \varepsilon)^{-2H})(1 - (t \vee r' \vee \varepsilon)^{-2H})$$
$$\times (1 - (s' \vee u \vee \varepsilon)^{-2H})(1 - (t' \vee u' \vee \varepsilon)^{-2H})$$
$$\times |r-r'|^{2H-2}|u-u'|^{2H-2}|s-s'|^{2H-2}$$
$$\times |t-t'|^{2H-2}\,dr\,dr'\,du\,du'\,ds\,ds'\,dt\,dt',$$

and one can show that the latter quantity is asymptotic to $\log(1/\varepsilon)$ for $\varepsilon$ converging to zero. This concludes the proof. $\square$

REMARK 4. An analog of Proposition 7 for the case $H = \frac{1}{2}$ is proved in Peccati and Yor (2004b).

3.2. *Quadratic functionals of a Brownian sheet.* We now extend the results of Deheuvels and Martynov (2004) and Peccati and Yor (2004a) to the case of a *Brownian sheet* $\mathbf{W}$ on $[0,1]^n$, that is,

$$\mathbf{W} = \{\mathbf{W}(x_1,\ldots,x_n) : (x_1,\ldots,x_n) \in [0,1]^n\}$$

is a centered Gaussian process with covariance function

$$\mathbb{E}[\mathbf{W}(x_1,\ldots,x_n)\mathbf{W}(y_1,\ldots,y_n)] = \prod_{i=1}^n (y_i \wedge x_i).$$

Note that the Gaussian space generated by $\mathbf{W}$ can be identified with an isonormal Gaussian process on $L^2([0,1]^n, dx_1\cdots dx_n)$.

In particular, we are interested in the limiting behavior of the two functionals, defined, respectively, for vectors $\beta = (\beta_1, \ldots, \beta_n)$ such that $\beta_i > -1$ and for $\varepsilon > 0$,

$$A_\beta = \int_{[0,1]^n} dx_1 \cdots dx_n \left( \prod_{i=1}^n x_i^{2\beta_i} \right) \mathbf{W}(x_1, \ldots, x_n)^2,$$

$$B_\varepsilon = \int_{[\varepsilon,1]^n} \frac{dx_1 \cdots dx_n}{(x_1 \cdots x_n)^2} \mathbf{W}(x_1, \ldots, x_n)^2,$$

when $\beta$ converges to $(-1, \ldots, -1)$ and $\varepsilon$ converges to zero. We recall that, due again to Jeulin's lemma, with probability 1,

$$\lim_{\beta \to (-1,\ldots,-1)} A_\beta = \lim_{\varepsilon \to 0} B_\varepsilon = +\infty.$$

REMARK 5. The law of the random variable $A_\beta$, for a fixed vector $\beta$ and for $n = 2$, is studied in detail in Deheuvels, Peccati and Yor (2004).

Now, for every $\beta$ and for every $\varepsilon$, we write the Wiener chaos decompositions of $A_\beta$ and $B_\varepsilon$ which are given by

$$A_\beta = \prod_{i=1}^n (2\beta_i + 2)^{-1} + \prod_{i=1}^n \frac{1}{2\beta_i + 1} I_2 \left( \prod_{i=1}^n (1 - (x_i \vee y_i)^{2\beta_i+1}) \right),$$

$$B_\varepsilon = \left( \log \frac{1}{\varepsilon} \right)^n + I_2 \left( \prod_{i=1}^n ((x_i \vee y_i \vee \varepsilon)^{-1} - 1) \right),$$

where $I_2$ now stands for a double Wiener integral with respect to $\mathbf{W}$. Then, calculations analogous to those performed in the previous section yield, thanks to Theorem 1, another CLT.

PROPOSITION 8. *Let the above notation and assumptions prevail.*

(i) *When $\beta \to (-1, \ldots, -1)$,*

$$\prod_{i=1}^n (2\beta_i + 2)^{1/2} \left( A_\beta - \prod_{i=1}^n (2\beta_i + 2)^{-1} \right)$$

$$= \prod_{i=1}^n (2\beta_i + 2)^{-1/2} \left( \prod_{i=1}^n (2\beta_i + 2) A_\beta - 1 \right)$$

*converges in distribution to $\sqrt{2^n} N(0,1)$, where $N(0,1)$ indicates a standard Gaussian random variable, and also $(\prod_{i=1}^n (2\beta_i + 2) A_\beta - 1)$ converges to zero in $L^2$.*

(ii) *When $\varepsilon \to 0$,*

$$\left(\log \frac{1}{\varepsilon}\right)^{-n/2} \left(B_\varepsilon - \left(\log \frac{1}{\varepsilon}\right)^n\right)$$

*converges in distribution to $\sqrt{2^n} N(0,1)$.*

**Acknowledgments.** Part of this paper was written while Giovanni Peccati was visiting the Institut de Matemàtica of the University of Barcelona (IMUB), in February 2003. This author heartily thanks David Nualart and Marta Sanz-Solé for their support and hospitality.

Facultat de Matemàtiques  
Universitat de Barcelona  
Gran Via 585  
08007 Barcelona  
Spain  
e-mail: nualart@mat.ub.es

Laboratoire de Statistique  
Théorique et Appliquée  
Université de Paris VI  
175 rue du Chevaleret  
75013 Paris  
France  
e-mail: giovanni.peccati@libero.it